\documentclass[12pt,a4paper]{amsart}
\usepackage{amsfonts, amsmath, amssymb, amsthm, hyperref}
\usepackage{anysize}

\newtheorem{thm}{Theorem}[section]
\newtheorem*{thm*}{Theorem}
\newtheorem{lem}[thm]{Lemma}

\theoremstyle{definition}
\newtheorem{defn}[thm]{Definition}

\theoremstyle{remark}

\numberwithin{equation}{section}

\DeclareMathOperator{\hind}{ind}
\DeclareMathOperator{\conv}{conv}

\DeclareMathOperator{\hdim}{hdim}

\newcommand{\pt}{{\rm pt}}
\renewcommand{\int}{\mathop{\rm int}}

\renewcommand{\epsilon}{\varepsilon}
\renewcommand{\phi}{\varphi}

\begin{document}
\title{A topological central point theorem}
\author{R.N.~Karasev}
\thanks{Supported by the Dynasty Foundation, the President's of Russian Federation grant MK-113.2010.1, the Russian Foundation for Basic Research grants 10-01-00096 and 10-01-00139, the Federal Program ``Scientific and scientific-pedagogical staff of innovative Russia'' 2009--2013, and the Russian government project 11.G34.31.0053.}

\email{r\_n\_karasev@mail.ru}

\address{Roman Karasev, Dept. of Mathematics, Moscow Institute of Physics
and Technology, Institutskiy per. 9, Dolgoprudny, Russia 141700}
\address{Roman Karasev, Laboratory of Discrete and Computational Geometry, Yaroslavl' State University, Sovetskaya st. 14, Yaroslavl', Russia 150000}


\begin{abstract}
In this paper a generalized topological central point theorem is proved for maps of a simplex to finite-dimensional metric spaces. Similar generalizations of the Tverberg theorem are considered.
\end{abstract}

\maketitle

\section{Introduction}

Let us state the discrete version of the Neumann--Rado theorem~\cite{neumann1945,rado1946,grun1960} (see also the reviews~\cite{eck1993} and~\cite{dgk1963}):

\begin{thm*}[The discrete central point theorem]
Suppose $X\subset\mathbb R^d$ is a finite set with $|X|=(d+1)(r-1)+1$. Then there exists $x\in\mathbb R^d$ such that for any halfspace $H\ni x$
$$
|H\cap X| \ge r.
$$
\end{thm*}

In this theorem a halfspace is a set $\{x\in \mathbb R^d : \lambda(x) \ge 0\}$ for a (possibly not homogeneous) linear function $\lambda : \mathbb R^d\to \mathbb R$. Using the Hahn--Banach theorem~\cite{rud1991} we restate the conclusion of this theorem as follows: the point $x$ is contained in the convex hull of any subset $F\subseteq X$ of at least $d(r-1)+1$ points.

When stated in terms of convex hulls, the central point theorem has an important and nontrivial generalization proved in~\cite{tver1966}:

\begin{thm*}[Tverberg's theorem]
Consider a finite set $X\in \mathbb R^d$ with $|X|=(d+1)(r-1)+1$. Then $X$ can be partitioned into $r$ subsets $X_1, \ldots, X_r$ so that
$$
\bigcap_{i=1}^r \conv X_i\neq\emptyset.
$$
\end{thm*}

In~\cite{barss1981,vol1996} a topological generalization of the Tverberg theorem was established. Instead of taking a finite point set in $\mathbb R^d$ and the convex hulls of its subsets, we take the continuous image of a simplex in $\mathbb R^d$ and the images of its faces (faces of the simplex viewed as a simplicial complex):

\begin{thm*}[The topological Tverberg theorem]
Let $m=(d+1)(r-1)$, $r$ be a prime power, and let $\Delta^m$ be the $m$-dimensional simplex. Suppose $f:\Delta^m\to Y$ is a continuous map to a $d$-dimensional manifold $Y$. Then there exist $r$ disjoint faces $F_1,\ldots, F_r\subset \Delta^m$ such that
$$
\bigcap_{i=1}^r f(F_i)\neq\emptyset.
$$
\end{thm*}

It is still unknown whether such a theorem holds for $r$ not equal to a prime power. But if we return to the central point theorem, we see that the following topological version holds without restrictions on $r$. Moreover, the target space can be any $d$-dimensional metric space, not necessarily a manifold. So the main result of this paper is:

\begin{thm}
\label{topcpt}
Let $m=(d+1)(r-1)$, let $\Delta^m$ be the $m$-dimensional simplex, and let $W$ be a $d$-dimensional metric space. Suppose $f:\Delta^m\to W$ is a continuous map. Then 
$$
\bigcap_{\substack{F\subset\Delta^m\\\dim F = d(r-1)}} f(F)\neq\emptyset,
$$
where the intersection is taken over all faces of dimension $d(r-1)$.
\end{thm}

Note that for $W=\mathbb R^d$ this theorem can also be deduced from the topological Tverberg theorem (see Section~\ref{proof-from-tv-sec} for details). The goal of this paper is to give another proof of Theorem~\ref{topcpt}, valid for any $d$-dimensional $W$. In Section~\ref{tver-sec} we show that a similar generalization of the Tverberg theorem for maps into finite-dimensional spaces essentially needs larger values of $m$.

\medskip
{\bf Acknowledgments. }
The author thanks Alexey~Volovikov, Pavle~Blagojevi\'c, Arseniy~Akopyan,  Peter~Landweber,  Vladimir Tikhomirov, and Satya~Deo for discussions and useful comments. 

\section{Index of $\mathbb Z_2$-spaces}

Let us recall some basic facts on the homological index of $\mathbb Z_2$-actions ($\mathbb Z_2$ is a group with two elements); the reader may consult the book~\cite{mat2003} for more details. Denote $G=\mathbb Z_2$, if we consider $\mathbb Z_2$ as a transformation group. The algebra $H^*(BG; \mathbb F_2)$ is a polynomial ring $\mathbb F_2[c]$ with the one-dimensional generator $c$. 

In this section we consider the cohomology with $\mathbb F_2$ coefficients, the coefficients being omitted from the notation. Define the equivariant cohomology for a space $X$ with continuous action of $G$ (a $G$-space) by
$$
H_G^*(X) = H^*(X\times_G EG) = H^*((X\times EG)/G),
$$
thus we have $H_G^*(\pt)=H^*(BG)$ for a one-point space with trivial action of $G$ and $H_G^*(X) = H^*(X/G)$ for a free $G$-space. For $G=\mathbb Z_2$ we may take $EG$ to be the infinite-dimensional sphere $S^\infty$ with the antipodal action of $G$, and $BG=\mathbb RP^\infty$. For any $G$-space $X$ the natural map $X\to \pt$ induces the natural cohomology map 
$$
\pi_X^* : H_G^*(\pt) = H^*(BG)\to H_G^*(X).
$$

\begin{defn}
For a $G$-space $X$ define $\hind_G X$ to be the maximal $n$ such that $\pi_X^*(c^n)\neq 0\in H_G^*(X)$.
\end{defn}

Note that if $X$ has a $G$-fixed point then the map $\pi_X^*$ is necessarily injective and the index is infinite. The following property of index is obvious by definition:

\begin{lem}
\label{ind-sup}
If $X$ is a topological disjoint union of $G$-invariant subspaces $X_1, \ldots, X_k$, then
$$
\hind_G X = \max_i \hind_G X_i.
$$
\end{lem}

The next property is the generalized Borsuk--Ulam theorem (see~\cite{mat2003} for example):

\begin{lem}
\label{odd-maps}
Let $\hind_G X\ge n$ and let $V$ be an $n$-dimensional vector space with antipodal $G$-action. Then for every continuous $G$-equivariant map $f:X\to V$
$$
f^{-1}(0)\neq \emptyset.
$$
\end{lem}

The following lemma is proved in~\cite{yang1955}, see also~\cite{kar2007}:

\begin{lem}
\label{even-maps}
Let $X$ be a compact metric $G$-space, $\hind_G X\ge (d+1)k$, and let $W$ be a $d$-dimensional metric space with trivial $G$-action. Then for every continuous $G$-equivariant map $f:X\to W$ there exists $x\in W$ such that
$$
\hind_G f^{-1}(x)\ge k.
$$
\end{lem}

In this lemma it is important to use the \v{C}ech cohomology, which is assumed in the sequel.

\section{Proof of Theorem~\ref{topcpt}}

Consider a continuous map $f:\Delta^m\to W$. Let us map the $m$-dimensional sphere $S^m$ to $\Delta^m$ by the formula:
$$
g(x_1,\ldots, x_{m+1}) = (x_1^2,\ldots, x_{m+1}^2).
$$
Apply Lemma~\ref{even-maps} to the composition $f\circ g$, which is possible because $g(x)=g(-x)$. We obtain a point $x\in W$ such that for $Z=(f\circ g)^{-1}(x)$ we have $\hind_G Z \ge r-1$.

We are going to show that $x$ is the required intersection point. Assume the contrary: a face $F\subseteq \Delta^m$ of dimension $d(r-1)$ does not intersect $g(Z)$. Without loss of generality, let $g^{-1}(F)$ be defined by the equations 
$$
x_1=\dots=x_{r-1} = 0.
$$
Note that the $r-1$ coordinates $x_1,\ldots, x_{r-1}$ give a continuous $G$-equivariant map $h : S^m\to \mathbb R^{r-1}$, where $G$ acts on $\mathbb R^{r-1}$ antipodally. By Lemma~\ref{odd-maps} the intersection $g^{-1}(F)\cap Z = h^{-1}(0)\cap Z = h|_Z^{-1}(0)$ should be nonempty. The proof is complete.

\section{Remark on the case $W=\mathbb R^d$ of Theorem~\ref{topcpt}}
\label{proof-from-tv-sec}

Recall the known fact: The case $W=\mathbb R^d$ of Theorem~\ref{topcpt} follows from the topological Tverberg theorem (only the case of prime $r$ is needed). For the reader's convenience we present a proof here (see also~\cite[Section~6]{kar2011}).

Consider a simplicial map $\phi : \Delta^M\to \Delta^m$, where $R=k(r-1)+1$ is a prime (for some $k$ this is so by the Dirichlet theorem on arithmetic progressions), $M=(R-1)(d+1)+k-1$, and there are $k$ vertices of $\Delta^M$ in the preimage of every vertex of $\Delta^m$. For $\Delta^M$ the topological Tverberg theorem holds (since $M\ge (R-1)(d+1)$), and there exist $R$ disjoint faces $\tilde F_1,\ldots, \tilde F_R$ of $\Delta^M$ such that
$$
\bigcap_{i=1}^R f(\phi(\tilde F_i))\ni x.
$$ 
Consider a face $F\subseteq \Delta^m$ of dimension $d(r-1)$ and assume that $\phi^{-1}(F)$ does not contain any $\tilde F_i$, then $M+1$ must be at least the number of vertices in $\phi^{-1}(F)$ plus $R$, that is
$$
M+1\ge k(r-1)d + k + R = (R-1)d+k+R = M+2,
$$
which is a contradiction. So $\phi^{-1}(F)$ contains some $\tilde F_i$, and $f(F)\ni x$.

\section{Tverberg type theorems for maps to finite-dimensional spaces}
\label{tver-sec}

It is natural to ask whether the corresponding version of the Tverberg theorem holds for maps from $\Delta^m$ to a $d$-dimensional metric space, at least for $r$ a prime power. In fact, the number $m=(d+1)(r-1)$ must be increased, as claimed by the following:

\begin{thm}
\label{tv-counter}
Let $m = (d+1)r-2$. Then there exists a $d$-dimensional polyhedron $W$ and a continuous map $f : \Delta^m\to W$ with the following property. For any pairwise disjoint faces $F_1,\ldots, F_r\subseteq \Delta^m$ there exists $i$ such that
$$
f(F_i)\cap f(F_j) = \emptyset
$$
for all $j\neq i$.
\end{thm}

This theorem also shows that our approach used to prove Theorem~\ref{topcpt} cannot be applied to the topological Tverberg theorem. Indeed, this proof does not distinguish between $\mathbb R^d$ and any metric $d$-dimensional space, but the topological Tverberg theorem does not hold for maps to $d$-dimensional metric spaces.

\begin{proof}[Proof of Theorem~\ref{tv-counter}]
The construction in the proof is taken from~\cite{volsce2005}. Let $\Delta^m$ be a regular simplex in $\mathbb R^m$, centered at the origin. Denote by $\Delta_{d-1}^m$ its $(d-1)$-skeleton, and $W = C\Delta_{d-1}^m$ the cone (geometrically centered at the origin) on this skeleton. Define the PL-map (of the barycentric subdivision to the barycentric subdivision) $f : \Delta^m \to W$ as follows. For every face $F\subseteq \Delta^m$ of dimension $\le d-1$ its barycenter is mapped to itself, for every face $F\subseteq \Delta^m$ of dimension $\ge d$ its barycenter is mapped to the origin.

Let $F_1,\ldots, F_r\subseteq \Delta^m$ be a set of $r$ pairwise disjoint faces. For some $i$ the dimension $\dim F_i$ is at most $d-1$ by the pigeonhole principle. For such a face we have $f(F_i)=F_i$, and
$$
f(F_i)\cap f(F_j)\subseteq F_i\cap f(F_j) \subseteq \partial\Delta^m.
$$
Since $f(F_j)\cap \partial\Delta^m \subseteq F_j$ we obtain
$$
f(F_i)\cap f(F_j)\subseteq F_i\cap F_j=\emptyset.
$$
\end{proof}

The following positive result for larger $m$ is a direct consequence of the reasoning in~\cite{vol2002}:

\begin{thm}
\label{tv-polyh}
Let $m=(d+1)r-1$ and let $r$ be a prime power. Suppose $f:\Delta^m\to W$ is a continuous map to a $d$-dimensional metric space $W$. Then there exist $r$ disjoint faces $F_1,\ldots, F_r\subset \Delta^m$ such that
$$
\bigcap_{i=1}^r f(F_i)\neq\emptyset.
$$
\end{thm}

\begin{proof}
Without loss of generality we may assume $W$ to be a finite $d$-dimensional polyhedron. Assume the contrary and denote $\Delta^m$ by $K$ for brevity. Then there exists a map 
$$
\tilde f : K^{*r}_{\Delta(2)} \to W^{*r}_{\Delta(r)}
$$
from the $r$-fold pairwise deleted join $K^{*r}_{\Delta(2)}$ in the simplicial sense to the $r$-fold $r$-wise deleted join $W^{*r}_{\Delta(r)}$ in the topological sense (see the definitions of the deleted joins in~\cite{mat2003}). Following~\cite{vol1996}, put $r=p^\alpha$ and consider the group $G=(\mathbb Z_p)^\alpha$ and let $G$ act on the factors of the deleted join transitively. The rest of the reasoning is based on the following facts from~\cite{vol2000,vol2002}:

Let $X$ be a connected $G$-space. Consider the Leray--Serre spectral sequence with
$$
E_2^{*,*} = H^*(BG; H^*(X; \mathbb F_p))
$$
converging to $H_G^*(X; \mathbb F_p)$. Here $G$ may act on $H^*(X; \mathbb F_p)$ so the cohomology $H^*(BG; \cdot)$ may be with twisted coefficients.

\begin{defn}
Denote by $i_G(X)$ the minimum $r$ such that the differential $d_r$ of this spectral sequence has nontrivial image in the bottom row.
\end{defn}

The index $i_G$ has the following properties, if $G$ is a $p$-torus $G=(\mathbb Z_p)^\alpha$:

\begin{enumerate}
\item (Monotonicity)
\label{monot}
If there is a $G$-map $f: X\to Y$, then $i_G(X) \le i_G(Y)$. If in addition $i_G(X)=i_G(Y)=n+1$ then the map $f^*:H^n(Y; \mathbb F_p)\to H^n(X; \mathbb F_p)$ is nontrivial.

\item (Dimension upper bound)
$i_G(X)\le \hdim_{\mathbb F_p} X+1$.

\item (Cohomology lower bound)
If $X$ is connected and acyclic over $\mathbb F_p$ in degrees $\le N-1$, then $i_G(X) \ge N+1$.
\end{enumerate}

Now note that from the cohomology lower bound it follows that $i_G\left(K^{*r}_{\Delta(2)}\right) \ge m+1$, from the dimension upper bound it follows that $i_G\left(W^{*r}_{\Delta(r)}\right)\le m+1$, and from (\ref{monot}) the map 
$$
\tilde f^* : H^m\left(W^{*r}_{\Delta(r)}; \mathbb F_p\right) \to H^m\left(K^{*r}_{\Delta(2)}; \mathbb F_p\right)
$$
must be nontrivial. From the cohomology exact sequence of a pair it follows that the natural map
$$
g^* : H^m\left(W^{*r}; \mathbb F_p\right) \to H^m\left(W^{*r}_{\Delta(r)}; \mathbb F_p\right)
$$
is surjective because $H^{m+1}(W^{*r}, W^{*r}_{\Delta(r)}; \mathbb F_p) = 0$ by dimensional considerations. Now it follows that the map
$$
(g\circ \tilde f)^* : H^m\left(W^{*r}; \mathbb F_p\right) \to H^m\left(K^{*r}_{\Delta(2)}; \mathbb F_p\right)
$$
is nontrivial. But the map $g\circ\tilde f$ is a composition of the natural inclusion
$$
h : K^{*r}_{\Delta(2)}\to K^{*r}
$$
with the map
$$
f^{*r} : K^{*r}\to W^{*r}.
$$
The latter map has contractible domain, and therefore induces a zero map on cohomology $H^m(\cdot; \mathbb F_p)$. We obtain a contradiction.
\end{proof}

\section{The case $r=2$ of Theorem~\ref{topcpt} and the Alexandrov width}

Let us give a definition, generalizing the definition in~\cite{tikh1976}. The reader may also consult the book~\cite{pink1985} in English. Throughout this section we use the notation 
$$
\delta A = \{\delta a: a\in A\}\quad\text{and}\quad A+B = \{a+b: a\in A,\ b\in B\}.
$$

\begin{defn}
\label{b-width-def}
Let $K\subseteq \mathbb R^n$ be a convex body. Denote by $b_k(K)$ the maximal number such that for any map $K\to Y$ to a $k$-dimensional polyhedron there exists $y\in Y$ such that for any $\delta < b_k(K)$ the set $f^{-1}(y)$ cannot be covered by a translate of $\delta K$.
\end{defn}

We use $k$-dimensional polyhedra $Y$ following~\cite{tikh1976}, but we may also use $k$-dimensional metric spaces as above. 

The definition of the \emph{Alexandrov width} (in~\cite{tikh1976}) is a bit different: A subset $X$ of some normed space $E$ is considered and $a_k(X)$ denotes the maximal number such that for any map $X\to Y$ to a $k$-dimensional polyhedron there exists $y\in Y$ such that for any $\delta < a_k(X)$ the set $f^{-1}(y)$ cannot be covered by a \emph{ball} (in the given norm of $E$) of radius $\delta$.

In~\cite[Theorem~1, p.~268]{tikh1976} the results of K.~Sitnikov and A.M.~Abramov~\cite{abr1972,sit1958} are cited, which assert that $a_k(X) = 1$ for any $k\le n-1$, if $X$ is the unit ball of a norm in $\mathbb R^n$. In terms of Definition~\ref{b-width-def} this means that $b_k(K) = 1$ for centrally symmetric convex bodies in $\mathbb R^n$ if $k \le n-1$ and obviously $b_k(K)=0$ for $k\ge n$.

Note that Theorem~\ref{topcpt} for $r=2$ actually asserts that $b_k(\Delta^n) = 1$ if $k\le n-1$. Indeed, if $f^{-1}(y)$ intersects all facets of $\Delta^n$ then it cannot be contained in a smaller homothetic copy of $\Delta^n$. Now it makes sense to extend the result of K.~Sitnikov and A.M.~Abramov to (possibly not symmetric) convex bodies:

\begin{thm}
\label{alex-waist}
If $K$ is a convex body in $\mathbb R^n$ and $k\le n-1$, then $b_k(K) = 1$.
\end{thm}

\begin{proof}
The proof in~\cite[Proposition~1, pp.~84--85]{tikh1976} actually works in this case. Assume the contrary: the map $f : K\to Y$ is such that every preimage $f^{-1}(y)$ can be covered by a smaller copy of $K$ and $\dim Y\le n-1$. For a fine enough finite closed covering of $Y$ its pullback covering $\mathcal U$ of $K$ has the following properties: the multiplicity of $\mathcal U$ is at most $n$ and any $X\in\mathcal U$ can be covered by a translate of $\delta K$ for some fixed $0<\delta <1$.

Assume $0\in \int K$ and call the point $t$ \emph{the center} of a translate $\delta K + t$. Assign to any $X\in\mathcal U$ the center $t_X$ of $\delta K + t_X \subseteq X$. Using the partition of unity subordinate to $\mathcal U$ we map $K$ to the nerve of $\mathcal U$, and then map this nerve to at most $(n-1)$-dimensional subcomplex of $\mathbb R^n$ by assigning $t_X$ to $X$. Finally we obtain a continuous map $\phi: K\to\mathbb R^n$ such that for any $x\in K$ we have $x\in \phi(x) + \delta K$ and the image $\phi(K)$ has dimension $\le n-1$.

Under the above condition the image $\phi(\partial K)$ cannot intersect $\varepsilon K$ if $\varepsilon < 1 - \delta$, because $\varepsilon K + \delta K$ is in the interior of $K$. If we compose $\phi|_{\partial K}$ with the central projection of $K\setminus \{0\}$ onto $\partial K$, we obtain a map homotopic to the identity map of $\partial K$. Therefore the map of pairs $\phi: (K, \partial K)\to (K, K\setminus \varepsilon K)$ has degree $1$, and $\phi(K) \supseteq \varepsilon K$. Therefore $\phi(K)$ is $n$-dimensional, which is a contradiction. 
\end{proof}

\end{document}